\documentclass[english]{article}

\usepackage{tikz}
\usepackage{pict2e}
\usepackage{pgfplots}
\usepackage{float}
\usepackage{url}
\usepackage{amsmath} 

\restylefloat{table}

\font\tenbb=msbm10
\font\sevenbb=msbm7
\font\fivebb=msbm5
\newfam\bbfam
\textfont\bbfam=\tenbb
\scriptfont\bbfam=\sevenbb
\scriptscriptfont\bbfam=\fivebb
\def\Bbb{\fam\bbfam\tenbb}

\begin{document}

\newtheorem{theo}{Theorem}
\newtheorem{prop}{Proposition}
\newtheorem{property}{Property}
\newtheorem{lemma}{Lemma}
\newtheorem{cor}{Corollary}
\newtheorem{defi}{Definition}
\newtheorem{rque}{Remark}

\title{On the size of the fibers of a morphism of algebraic varieties
}
\author{Abdelkader Mokkadem\footnote{Universit\'e Paris-Saclay, UVSQ, CNRS, Laboratoire de Math\'ematiques de Versailles, 78000, Versailles, France (abdelkader.mokkadem@uvsq.fr)}}
\date{}
\maketitle

\begin{abstract}
 Let $ u : V \rightarrow W $ be a morphism of algebraic  varieties. Given an integer $m$ or $m = \infty$, we prove that the set of $y \in W$ such that $u^{-1}(y)$ has size $m$, is a constructible set.  We also show that this set is locally closed when $W$ is normal and $u$ is finite.

\end{abstract}

.




\section{Introduction}

Let $ u : V \rightarrow W $ be a morphism of algebraic  varieties.  We will denote by $\vert A \vert $ the size of a set $A$.  The size is either an integer or infinite. For each $y \in W $ the fiber is $ u^{-1}(y)$ . The goal of this paper is to study how the size of the fibers is distributed in W. For this, we will focus in particular on the sets  $ \{ y \in W, \  \  \vert u^{-1}(y) \vert = m \} $ where $m \in {\Bbb N}$ or $m = \infty $. In the following, the set $ \{ y \in W, \  \  \vert u^{-1}(y) \vert = m \} $  will be denoted by $F_{m} $. One of the questions we address is whether $F_{m} $ is constructible.

For the definition and properties of  constructible set, we refer to \cite{Humphreys} and \cite{Mumford}. Briefly, a set is said to be constructible if it is a finite union of locally closed sets—a locally closed set being the intersection of an open set and a closed set.

The definitions and properties of affine or non affine algebraic varieties over a field $k$ are given in \cite{Görtz}, \cite{Hartshorne}, \cite{Humphreys}, \cite{Mumford} and \cite{Perrin}. Throughout this paper $k$ is an algebraically closed field  and we will specify each time when the variety is irreducible. A non-irreducible (affine or not) algebraic variety $Z$ is simply a union of finitely many irreducible (affine or not) algebraic varieties. In this case $Z$ can be written uniquely as an irredundant union of finitely many irreducible algebraic varieties called the components of $Z$. Of course, if $Z$ is affine the components are affine.  We will always consider an algebraic variety as a topological space and the only topology we will use will be the Zariski topology.

The distribution of fibers according to a characteristic is an important element in the study of morphisms. 

Some characteristics that are often studied are properties. More precisely, given a property (P), we study whether  the fiber satisfies the property (P) or not. The study of the set of all $y \in W $ such that the fiber of $y$ possesses a property (P) is considered within the framework of Grothendieck’s localization problem (see \cite{Grothendieck}). For references and  results on this subject see \cite{Avramov} and \cite{Shimomoto}.

The most commonly used characteristic of the same type as size, is dimension. In this case, we have the famous  Semi-Continuity Theorem of Chevalley (see \cite{Humphreys} or  \cite{Mumford}). 
This theorem states that for any integer $m$, the set  $ \{ x \in V, \  \ dim( u^{-1}(u(x)) ) \geq m  \} $ is a closed subset of $V$. It follows from this theorem that for any integer $m$, the set $ \{ x \in V, \  \ dim( u^{-1}(u(x)) ) = m  \} $ is a constructible subset of V. 

Now, we know that the image of a constructible set under the morphism $u$ is a constructible set (by virtue of another important theorem, the Chevalley's theorem on constructible sets, see \cite{Humphreys}  or  \cite{Mumford} ). We deduce that for any integer $m$, 
the set $ \{ y \in W, \  \ dim( u^{-1}(y) ) = m \} $ is a constructible subset of $W$. 

In our study of the size of the fibers, we obtain a similar result. We consider special cases as well as the most general case. In the very important special case where $V$ and $W$ are irreducible, $W$ is normal and $u$ is dominant and finite with separable degree $n$, the result is very simple and explicit. We will show that for any integer $ s $, the set  $ \{ y \in W, \  \ \vert u^{-1}(y) \vert \geq s \} $ is an open subset of $W$ which we define precisely. This set is empty for $s > n$. 
This implies in particular that $F_{n}$ is open and  that for any integer $ s $, the set  $F_{s} $ is a locally closed set. 

In the general case we are considering, $V$ and $W$ are algebraic varieties that are assumed to be neither affine nor irreducible. We then show that, for any integer $s$ or for $s = \infty$, the set $F_s$ is constructible.

There are no results on the size of fibers in the published literature. The only place where we have found the issue of fiber size addressed is in  Milne’s Algebraic Geometry notes course (\cite{Milne}). In page 184, he provides results on the size in the case where $u$ is finite surjective and $W$ is normal.

Our paper is organized as follows. Section 2 is devoted to some preliminaries. Section 3 deals with the case where $W$ is normal and $u$ is finite. Section 4 addresses the general case.

\section{Preliminaries}

\subsection{The irreducible and dominant case} \label{Decomposition}

The irreducible and dominant case is essential in all the results. Let us present the framework in this case.

Recall that throughout this article, $k$ denotes an algebraically closed field.  Let $V$ and $W$ be irreducible affine algebraic  varieties over $k$.  We denote by  $k[V]$ the ring of regular functions on $V$.  Since $V$ is irreducible $k[V]$ is an integral domain. We denote its field of fractions by $k(V)$. In the same way, we denote  by $k[W]$  the ring of regular functions on $W$ and by $k(W)$  the field of fractions of  $k[W]$.
Let $ u : V \rightarrow W $ be a dominant morphism, the associated k-algebra homomorphism  $ \varphi : k[W] \rightarrow k[V]$, $ \varphi (f) = f\circ u $, is injective. The homomorphism $ \varphi$ makes  $k[V] $  a $k[W]$-algebra of finite type, $k[V] = k[W][H_{1},\ldots,H_{n}]$.  Since  $ \varphi$ is injective,   we extend $ \varphi$ to a map on the fields $k(W)\rightarrow k(V)$. When $r= dim(V)- dim(W)=0 $, $k(V)$ is algebraic over $k(W)$, $[k(V): k(W)]$ is the degree of $u$ and $[k(V): k(W)]_{s}$ is the separable degree of $u$.  
When $r > 0 $, the transcendence degree of  $k(V)$  over $k(W)$ is $r$.

Regarding the fibers of $u$, it is known (see \cite{Perrin}) that for all $y \in u(V)$ and all irreducible component $Z$ of $u^{-1}(y)$, $\dim(Z) \geq r$.  Moreover there exists a non-empty open subset $U$ of $W$ such that, for all $y \in U$, $\dim(u^{-1}(y)) = r$. 
In particular, in the case where $r > 0$, all fibers are either empty or infinite. 

To study the morphism $u$, we will often  write it as a composition of morphisms. Finding a decomposition of $u$ is equivalent to finding a decomposition of the injection $ \varphi : k[W] \rightarrow k[V] $.  A  way to find decomposition is as follows. Let $H_{1},\ldots,H_{n}$ be a  set of generators of $k[V]$ over $ k[W] $,  $k[V] = k[W][H_{1},\ldots,H_{n}]$. Let  $ \{ H_{i_{1}},\ldots,H_{i_{s}} \}$ a subset of $ \{ H_{1},\ldots,H_{n} \}$. We then have   $ \varphi = \nu \circ \psi $ where $ \psi :  k[W] \rightarrow k[W][H_{i_{1}},\ldots,H_{i_{s}}]$ is the injection defined by $ \psi (z) = \varphi (z) $ and $ \nu : k[W][H_{i_{1}},\ldots,H_{i_{s}}] \rightarrow k[V] $ is the natural injection $ \nu (z) = z $. 
This decomposition of $ \varphi $ corresponds to a decomposition of $u$. The domain $k[W][H_{i_{1}},\ldots,H_{i_{s}}]$ is the ring of regular functions of an irreducible algebraic affine variety $V_{0}$. The injection $ \psi $ corresponds to a dominant morphism $ u_{0} : V_{0}  \rightarrow W $, the injection $ \nu $ corresponds to a dominant morphism $ \mu : V  \rightarrow V_{0} $ and we have  $ u =  u_{0} \circ \mu$. 
Since the morphisms are dominant, we have $dim(W) \leq dim(V_{0}) \leq dim(V) $.

\subsection{Injective morphism}

We present here a simple but very useful result. 

\begin{prop} [The injective morphism] \label{Injective morphism} 
Let $V$ and $W$ be irreducible affine algebraic  varieties over $k$ and let $ u : V \rightarrow W $ be a dominant morphism. 
Let $H_{1},\ldots,H_{n}$ be such that $k[V] = k[W][H_{1},\ldots,H_{n}]$ and let  $ \Psi = (u, H_{1},\ldots,H_{n}) $. Then $\Psi$ is a regular injective morphism from $V$ to $W \times k^{n} $. 

\end{prop}

\paragraph{Proof}
Of course, if $n=0$, $ \Psi = u $ is an isomorphism. We consider the case $n\geq1$.

For every $F \in k[V] = k[W][H_{1},\ldots,H_{n}] $  , there are $C_{1},\ldots, C_{N}$ in $ k[W]$ such that for every $x \in V$  
$$ F(x) = C_{1}(u(x))P_{1}(x) + \ldots+ C_{N}(u(x))P_{N}(x) $$
where $ P_{k} = H_{1}^{i_{1k}}\ldots H_{n}^{i_{nk}}$ and  $ i_{1k},\ldots,i_{nk} $ are non-negative integers.

Now let $x,x'$ in $V$, such that $u(x)=u(x'),  \ and \  H_{i}(x)=H_{i}(x'), i=1,\ldots,n$.  It follows that for every $F \in k[V] $, $F(x)=F(x')$. 
This is not possible if $ x \neq x'$ since for all  $x \neq x'$, there exists at least one coordinate function $X_{i} \in  k[V] $ such that $X_{i}(x) \neq X_{i}(x')$.

\paragraph{Remark}

Let $V'$  the closure of $ \Psi (V) $ in $W \times k^{n} $. Since $ \Psi $ is injective, we have $ dim(V') = dim(V). $

\subsection{Reduction to the separable case}
 For the definitions and properties of inseparable and purely inseparable extensions, we refer to \cite{Lang}. 
 The separable degree of $E$ over $F$ is noted $[E;F]_{s}$.
 The following proposition reduces the inseparable case to the separable case.

\begin{prop} [The reduction] \label{reduction}
Let $V$ and $W$ be irreducible affine algebraic  varieties over $k$ such that $ dim(V) = dim(W) $  and let $ u : V \rightarrow W $ be a dominant morphism.  
There exist an  irreducible affine algebraic  variety $\hat{V}$ and dominant morphisms $\hat{u}_{0} $ and $ \hat{u}_{1} $, $ \hat{u}_{0} :  V \rightarrow \hat{V} $, $ \hat{u}_{1} :  \hat{V} \rightarrow W $ such that 
\begin{itemize}
\item[$(a)$]
 $dim(\hat{V})= dim(W) $, $ u = \hat{u}_{1} \circ \hat{u}_{0} $, $k(\hat{V})$  is separable over $k(W)$ and  $[k(\hat{V}):k(W)]=[k(V): k(W)]_{s}$
\item[$(b)$] $\hat{u}_{0} $ is finite and bijective
\item[$(c)$]  for any $ y \in W $ we have $ \vert u^{-1}(y) \vert = \vert \hat{u}_{1}^{-1}(y) \vert $
\item[$(d)$]   $u$ is finite  if and only if  $ \hat{u}_{1}$ is finite.
\end{itemize}
\end{prop}

\paragraph{Proof}

Since $ dim(V) = dim(W) $, $k(V)$ is an algebraic extension over $k(W)$.  There are $H_{1},\ldots,H_{n}$, such that  $k[V] = k[W][H_{1},\ldots,H_{n}]$. Let $p$ be the characteristic of $k$. If the extension is separable, we have $\hat{V}=V$, $\hat{u}_{0}=Id $ and $ \hat{u}_{1}= u$. Let us assume the extension is not separable and therefore $p>0$. 

$(a)$  For each $i$, there exists an integer
 $ \mu_{i} \geq  0 $ such that $ H_{i}^{p^{\mu_{i}}}$ is separable over $k(W)$ (see \cite{Lang}). Let $i_{1},\ldots,i_{l}$ be the indices $i$ such that $ \mu_{i} > 0 $.  Let $\hat{V}$ the irreducible affine algebraic  variety associated to the integral domain $k[\hat{V}] = k[W][H_{1}^{p^{\mu_{1}}},\ldots,H_{n}^{p^{\mu_{n}}}]$. Clearly $ dim(\hat{V}) = dim(W) $. We decompose $u $ by using  the  injections
$$  k[W] \hookrightarrow k[\hat{V}] \hookrightarrow k[V]. $$
 We obtain $ u = \hat{u}_{1} \circ \hat{u}_{0} $ where $\hat{u}_{0} $ and $ \hat{u}_{1} $ are dominant.  For the morphism  
$ \hat{u}_{1} :  \hat{V} \rightarrow W $, the field $k(\hat{V})$  is separable over $k(W)$ and the degree of the extension is the separable degree of $k(V)$ over $k(W)$.

$(b)$  For the morphism $ \hat{u}_{0} :  V \rightarrow \hat{V} $, $k(V)=k(\hat{V})[H_{i_{1}},\ldots,H_{i_{l}}]$ is purely inseparable over $k(\hat{V})$ (see \cite{Lang}). For each $ i_{j} \in \{  i_{1},\ldots,i_{l} \} $,  $H_{i_{j}}$ is the unique multiple root of its irreducible polynomial $P(T) = T^{p^{\mu_{i_{j}}}} - a$, $a \in k[\hat{V}]$, $a = H_{i_{j}}^{p^{\mu_{i_{j}}}}$.  It follows that  $ \hat{u}_{0}$ is finite and therefore surjective. Moreover, for each $y  \in \hat{V}$, $y= \hat{u}_{0}(x)$, there is a unique $ z \in k$ such that $ z^{p^{\mu_{i_{j}}}} - a (y) =0 $. Applying the  Proposition~\ref{Injective morphism}, we show that $ \hat{u}_{0}$ is injective and therefore bijective.

$(c)$  It is obvious that, for any $ y \in W $ we have $ \vert u^{-1}(y) \vert = \vert \hat{u}_{1}^{-1}(y) \vert $.

$(d)$ is obvious.

\paragraph{Notation}

In what follows, to any dominant morphism $ u : V \rightarrow W $, where $V$ and $W$ are irreducible affine varieties and $ dim(V) = dim(W) $, we will associate  $\hat{V}$,  $\hat{u}_{1}$ and  $\hat{u}_{0}$  obtained by Proposition \ref{reduction}. In the separable case, $\hat{V}=V$, $\hat{u}_{0}=Id $, $ \hat{u}_{1}= u$ and we will say that $u$ is separable.

\subsection{Roots of polynomial}  .

Let $W$ be an irreducible affine algebraic  variety over $k$ and let $F(T) \in k[W][T]$ a monic polynomial, $F(T)=a_{0}+ a_{1}T +\ldots+a_{n-1} T^{n-1}+T^{n} $, $ a_{j} \in k[W]$. Let $H_{1}, H_{2},\ldots,H_{n}$  be the roots of $F$ and let  $A = k[W][H_{1}, H_{2},\ldots,H_{n}]$ be the integral domain generated by $H_{1}, H_{2},\ldots,H_{n}$. Let $\mathcal{W}$ be the irreducible affine algebraic variety whose ring of regular functions is $A$. To the  injection $  k[W] \hookrightarrow A $ corresponds a dominant morphism $ \mu : \mathcal{W} \rightarrow W $. Since $A$ is integral over $  k[W]$, $ \mu $ is a finite morphism and thus closed and surjective. 

For $y \in W $, we  denote $ F(y, T) = a_{0}(y)+ a_{1}(y)T +\ldots+a_{n-1}(y) T^{n-1}+T^{n} $. The set of all roots of $ F(y, T)$ is denoted $ Root(F, y)$,  therefore $ Root(F, y) = \{ z \in k, \ \ F(y, z) =0 \} $. For any $ x \in \mathcal{W} $ and $H \in A $, we have $ F(\mu(x), H(x)) = a_{0}(\mu(x))+ a_{1}(\mu(x))H(x) +\ldots+a_{n-1}(\mu(x)) H^{n-1}(x)+H^{n}(x)$. The roots of $ F(\mu(x), T)$ are $H_{1}(x), H_{2}(x),\ldots,H_{n}(x)$. For $ m \geq 1 $, 
we set $\mathcal{W}_{m}(F) = \{ x \in \mathcal{W}, \ \vert Root(F, \mu(x)) \vert < m \}$. It is easy to see that $\mathcal{W}_{m}(F)$ is a closed  subset of  $\mathcal{W} $. In particular, for $m=1$ it is empty and for $m > n$ it is $\mathcal{W}$. Let us set $ W_{m}(F) = \{ y \in W, \ \vert Root(F, y) \vert < m \}$. Clearly $ \mu (\mathcal{W}_{m}(F))= W_{m}(F)$, $ \mu^{-1} (W_{m}(F)) =  \mathcal{W}_{m}(F) $ and $W_{m}(F)$ is a closed  subset of $W$. We denote by $ U_{m}(F,W)$ the  open set  $ W \setminus W_{m}(F)$. The sequence $ U_{m}(F,W)$ is decreasing that is $ U_{m+1}(F,W) \subset  U_{m}(F,W)$. Of course $ U_{m}(F,W)$ is empty for $m > n $  and $ U_{1}(F,W) = W$.

\section{Normal finite morphism}

We start with the normal case, because in this case, the distribution of sizes is very  explicit. In particular the function $ size(y) = \vert u^{-1}(y) \vert $,  defined on W, is lower semi-continuous.

Let us specify the setting. Let $V$ and $W$ be irreducible affine algebraic varieties over $k$.  
We consider the case where  $W$ is normal and $ u : V \rightarrow W $ is a finite dominant morphism. It follows that  $ dim(V) = dim(W) $, that $u$ is closed and surjective and that for all $ y \in W$, $ 1 \leq \vert u^{-1}(y) \vert < \infty $. In particular, $F_{0}$ and  $F_{\infty}$ are empty.

Let $ n = [k(V):k(W)]_{s} $.  
 Let $\hat{V}$ be given by Proposition \ref{reduction} and let   $ \mathcal{P} $ be the set of monic polynomials  $P$ with degree $n$, such that $P$ is the minimal polynomial of some $\xi  \in k[\hat{V}] $ over $k(W)$. Clearly we have  $ k(\hat{V})= k(W)[\xi] $ since the degree of the extension $ k(\hat{V})$ is $n$   and  $ P \in k[W][T] $ since $W$ is normal.
 
 For all $ m \geq 1 $,   let $ U_{m}(u) = \cup_{P \in \mathcal{P}}  U_{m}(P,W)$. Of course by Noetherianity, there are $P_{1}$,\ldots,$P_{s}$ in  $\mathcal{P}$ such that   $ U_{m}(u) = U_{m}(P_{1},W) \cup \ldots\cup  U_{m}(P_{s},W) $. Since the irreducible minimal polynomials $P$ are separable, the  open subset $ U_{m}(u)$ is empty if and only if $m>n$. Moreover $ U_{1}(u)=W$.

\begin{theo} [Normal case] \label{Normal Case}
Let $W$ and  $V$ be irreducible affine algebraic  varieties and $ u : V \rightarrow W $  a finite dominant morphism. Let $ n = [k(V):k(W)]_{s} $ and assume that $W$ is normal. Then 
\begin{itemize}
\item[(i)] for all $ y \in W$, $ \vert u^{-1}(y) \vert \leq n $ 

\item[(ii)] $ \{ y \in W, \  \vert u^{-1}(y)\vert =n \} = U_{n}(u) $ 

\item[(iii)] for all  $ m \geq 1 $, 
$\{ y \in W, \    \vert  u^{-1}(y) \vert \geq m \} = U_{m}(u) $.
\end{itemize}

\end{theo}

\paragraph{Proof}  

 In the purely inseparable case,  $n=1$, $u$ is bijective, and $U_{m}(u)$ is empty for $m > 1$. The result is therefore trivial.

Let us now consider the general case. We note that (iii) implies (i) and (ii) since $ U_{n+1}(u)$ is empty. Only (iii) remains to be proved.   (iii) can be rewritten as $\{ y \in W, \ \   \vert  \hat{u}_{1}^{-1}(y) \vert \geq m \} = U_{m}(u)$.

\paragraph{Step 1}

We first prove that  $\vert \hat{u}_{1}^{-1}(y) \vert \geq m  \Rightarrow y \in U_{m}(u) $. 

The case $m=1$ is trivial since $ U_{1}(u)=W$. We consider the case $m >1$.
Let $ y \in W $ and $ \hat{u}_{1}^{-1}(y) = \{x_{1},\ldots,x_{l} \} $ where the $x_{j}$ are distinct and $ l \geq m $. 
 We can find $a \in k[\hat{V}] $ such that $ a(x_{i}) \neq  a(x_{j})$ if $ i \neq j $ . Here is one way to find such a regular function. We can view $\hat{V}$ as a subset of $ k^{d}$.  It therefore suffices to find a regular function $a$ on $ k^{d}$ that is injective on the finite subset $ \hat{u}_{1}^{-1}(y)$. In fact, we can choose a linear form for $a$. More precisely, the  $ x_{i} - x_{j } , i>j $, are non-zero in $ k^{d}$. We can then easily prove by induction on $d$ that there exists a linear form $a$ such that $ a( x_{i} - x_{j }) \neq 0 , i>j $. This yields the regular function that is injective on $ \hat{u}_{1}^{-1}(y)$.
 
We will now show that there exists $b \in k[\hat{V}]$ whose minimal polynomial has degree $n$ and such that $b$ is injective on the fiber $\hat{u}_{1}^{-1}(y)$.  Let $b_{0} \in k[\hat{V}] $  such that $ k(\hat{V})= k(W)[b_{0}] $. Since $k(\hat{V})$ is separable over $k(W)$, there exists only a finite number of fields $F$ such that $k(W) \subset F \subset k(\hat{V})$. It is then easy to show that there exists a finite subset $N$ of $k$, such that, for all $ \lambda \in k \setminus N $, we have $ k(\hat{V})= k(W)[a + \lambda b_{0}] $. 
On the other hand, let $M = \{ \dfrac{a(x_i) - a(x_j)}{b_{0}(x_j) - b_{0}(x_i)} , 1 \leq j < i \leq l, \, b_{0}(x_j) - b_{0}(x_i) \neq 0 \}$. $M$ is a finite subset of $k$  and, for all $ \lambda \in k \setminus M $, $b = a + \lambda b_{0}$ is injective on the fiber $\hat{u}_{1}^{-1}(y)$.

Consequently, we can choose $ \lambda \in k $ such that  $b = a + \lambda b_{0}$ is injective on the fiber $\hat{u}_{1}^{-1}(y)$ and has a minimal polynomial of degree $n$. 

Let $P(T)=a_{0}+ a_{1}T +\ldots+a_{n-1} T^{n-1}+T^{n} $ be the minimal polynomial of $b$. The values $b(x_{1}),\ldots,b(x_{l})$ are distinct roots of the polynomial $ P(y, T) = a_{0}(y)+ a_{1}(y)T +\ldots+a_{n-1}(y) T^{n-1}+T^{n} $. Since $l \geq m $ we therefore have $y \in U_{m}(P,W)$. 
 
 This proves that $\vert \hat{u}_{1}^{-1}(y)\vert \geq m  \Rightarrow y \in U_{m}(u) $.
 
 \paragraph{Step 2}

Let us now prove the converse.

Let $P \in  \mathcal{P} $  and $y \in U_{m}(P,W)$. We must show that $\vert \hat{u}_{1}^{-1}(y) \vert \geq m $.  Let $b \in k[\hat{V}] $  have minimal polynomial $P$ and let $V'$ be the irreducible affine algebraic variety whose ring of regular functions is $k[W][b]$. The injections $ k[W]\hookrightarrow k[W][b] \hookrightarrow k[\hat{V}] $ induce dominant morphisms $ \nu' : \hat{V} \rightarrow V' $ and $\nu : V' \rightarrow W $ such that $ \hat{u}_{1} =\nu \circ \nu' $. Since the morphism $ \nu' $ is finite, it is  surjective, and therefore $\vert \hat{u}_{1}^{-1}(y) \vert \geq \vert \nu^{-1}(y) \vert $. It therefore suffices to show that $\vert \nu^{-1}(y)\vert \geq m $.

 The map $ \varphi : V' \rightarrow W \times k $, $ \varphi(x) = (\nu(x),b(x)) $, is injective. Since, for any $x\in V' $, $b(x)$ is a root of $ P(\nu(x), T)$, we have $ \varphi (V') \subset W'$ where $W' = \{ (y,z) \in W \times k, \ \  P(y,z)=0 \} $ . It is clear that $W'$ is an irreducible affine variety and that $dim(W')=dim(W)=dim(\hat{V}) $. Consequently, the injective regular map $ g : V' \rightarrow W' $, where $g(x)= (\nu(x), b(x)) $, is dominant. Since $g$ is finite, it is surjective and therefore bijective. Let $ \pi : W' \rightarrow W $,  $ \pi (y,z) = y $. We have $\nu(x) =  \pi (g(x)) $ for every $x$ and, consequently,   
 $\nu^{-1}(y)= g^{-1}(\pi ^{-1} (y)) $. Since  $\pi ^{-1} (y) = \{(y,z),  z \in k,\  P(y,z)=0 \} $, we have   $   \vert \nu^{-1}(y) \vert \geq \vert \{ z \in k,\   P(y,z)=0 \} \vert $ and therefore $ \vert \nu^{-1}(y) \vert \geq m $.  

\paragraph{Comments}
$a)$ It follows from $(iii)$ that the function $ size(y) = \vert u^{-1}(y) \vert $ is lower semi-continuous.

$b)$ Part $(i)$ of the theorem says that $F_{m}$ is empty for $m > n$. Part $(ii)$ says that $F_{n} $ is a  non-empty open subset of $W$.\\

Let us specify the nature of all the $F_{m}$.  First, let us recall that a subset $C$ of $W$ is locally closed if $C = U \cap Z $ where $U$ is open and $Z$ is closed in $W$. Note that $C$ is locally closed if it is an open subset of a closed set. In particular, $C$ is locally closed if it is open within its closure. 

A constructible set is a finite union of locally closed sets. A locally closed set is therefore a particular constructible set. 

The following corollary is  obvious.

\begin{cor} [Normal Case] \label{Normal Case}
Under the assumptions of Theorem \ref{Normal Case}, for any  $ 1 \leq m \leq n  $, the set $F_{m} = \{ y \in W,\ \ \vert u^{-1}(y) \vert =m \}$ is a locally closed subset of~$W$. 

\end{cor}

\section{The general case}

We now consider the general case. We do not assume that the algebraic varieties $V$ or $W$ are affine or irreducible.

\begin{theo} [General Case] \label{General Case}
Let $V$ and $W$ be algebraic varieties and let  $ u : V \rightarrow W $ be a  regular morphism.
For any $m$ integer or $m = \infty$, $F_{m}$ is a constructible set.   Moreover, there is $l$ such that $F_{j}$  is empty for all $j$ such that $ l < j <  \infty $.

\end{theo} 

Let $W'$ be the closure of $u(V)$. The set $ W \setminus W' $ is constructible and $F_{0}$ is the disjoint union of  $ W \setminus W' $ and $ W' \cap  F_{0}$. Therefore $F_{0}$ is constructible if and only if $ W' \cap  F_{0}$ is constructible. Moreover, for $m \neq 0$, $F_{m} \subset W'$.
We can therefore assume that $u$ is dominant.

We first consider the set $F_{0}$.

\subsection{The set $F_{0}$}

Let us set $E_{0} = W \setminus F_{0}=u(V)$. 
We have

\begin{prop}  \label{ The image }

$F_{0}$ and $E_{0}$ are constructible.

\end{prop}

\paragraph{Proof}

Since $F_{0}$ and $E_{0}$ are complementary, if one of them is constructible, the other one is also constructible. 
The set $E_{0}$ is constructible by virtue of  Chevalley's theorem (see \cite{Humphreys} or \cite{Mumford}), and we can therefore deduce that $F_{0}$ is constructible. 
However, we will directly show that $F_{0}$ is constructible and thus deduce that $E_{0}$ is constructible. This provides a new proof of Chevalley's theorem.

To show that $F_{0}$ is constructible, we will use the following property.

\paragraph{(P)}

Let $Z$ be a closed subset of $W$ and let $C_{1},\ldots,C_{s}$ the irreducible components of $Z$. If  $Z \cap F_{0}$ is dense in $Z$ then for each $i$ the closure of $C_{i} \cap E_{0}$ is a proper subset of $C_{i}$. In particular the closure of $Z \cap E_{0}$ is a proper subset of $Z$ and $ dim (Z \cap E_{0}) < dim(Z) $.

Let us show the property \textbf{(P)}. Noting that  $Z \cap F_{0}$ is dense in $Z$ if and only if for each $i$, $C_{i} \cap F_{0}$ is dense in $C_{i}$, it suffices to show \textbf{(P)} in the case $Z$ is irreducible.  Let us therefore show \textbf{(P)} when $Z$ is irreducible.   

Let us assume the closure of $Z \cap E_{0}$ is not a proper subset of $Z$ and let us show that this implies that $Z \cap F_{0}$ is not dense in $Z$. Since in this case the closure of  $Z \cap E_{0}$ is $Z$, it follows that there is a component $L$ of $u^{-1}(Z)$ such that the restriction $ u : L \rightarrow Z$ is a dominant morphism. We then deduce that $u(L)$ contains a non-empty open subset of $Z$. Since $Z \cap F_{0}$ is contained in the complement of $u(L)$ in $Z$, it is therefore not dense in $Z$. This shows \textbf{(P)}.\\

Let us now complete the proof of the proposition. Let $Z_{0}$ be the closure of $F_{0}$. We have $ dim(Z_{0}) < dim (W)$. Since $Z_{0} \cap F_{0}$ is dense in $Z_{0}$, we apply \textbf{(P)}. 
Let $W_{1}$ be the closure of $ E_{0} \cap Z_{0}$. We have $ dim(W_{1}) < dim(Z_{0}) $ and $ F_{0}= \{Z_{0} \setminus W_{1} \} \cup F_{01}$ where $ F_{01} = \{ F_{0} \cap W_{1}\} $. 
We note that $ \{Z_{0} \setminus W_{1} \} $ is constructible.

Now, we proceed with $W_{1}$ and $ F_{01}  $ as we did with $W$ and $F_{0}$. 

Let $Z_{1}$ be the closure of $F_{01}$. Since $ E_{0} \cap Z_{0}$ is dense in $W_{1}$, we have $ dim(Z_{1}) < dim (W_{1})$. Now since $Z_{1} \cap F_{0}$ is dense in $Z_{1}$, we apply \textbf{(P)}.  Let $W_{2}$ be the closure of $ E_{0} \cap Z_{1}$. We have $ dim(W_{2}) < dim(Z_{1}) < dim (W_{1}) $ and $ F_{0}= \{Z_{0} \setminus  W_{1} \} \cup \{Z_{1} \setminus  W_{2} \}  \cup F_{02} $ where $ F_{02} = \{ F_{0} \cap W_{2}\} $.
We continue in this way $k$ times. We obtain $ F_{0}= \{Z_{0} \setminus  W_{1} \} \cup\ldots\cup  \{Z_{k-1} \setminus  W_{k} \}  \cup F_{0k} $  where $ F_{0k} = \{ F_{0} \cap W_{k}\} $. Since for all $i$, $ dim(W_{i+1}) < dim(W_{i}) $, there is $k$ such that $ W_{k} $ is empty and therefore 
$ F_{0}= \{Z_{0} \setminus  W_{1} \} \cup\ldots\cup  \{Z_{k-1} \setminus  W_{k} \}$.  Since each $ \{Z_{i-1} \setminus  W_{i} \} $ is constructible, $ F_{0}$ is also constructible. \\

To prove that $F_{m}$ is constructible when $m \neq 0$, we need some preliminary results. 

We first consider the case of a very specific finite morphism: the strong finite case.
We then deduce the result for an open subset $U$ of $W$:  the local result.
In a next step, we will see how the local result on $W$ extends to the closed subset $W \setminus U$: the result on closed subsets.  

\subsection{The case of strongly finite morphism}

Let $V$ and $W$ be irreducible affine algebraic  varieties over $k$  and let $ u : V \rightarrow W $ be a  dominant morphism. 

 We will say that the morphism $u$ is strongly finite  if there exists $\hat{V}$ given by Proposition \ref{reduction} such that   $k[\hat{V}] = k[W][H]$ where $H $ is  algebraic over $k(W)$ with minimal polynomial $G(T) \in k[W][T]$.   It follows that  $H$ is integral over $k[W]$, $u$ is finite and $ dim(V) = dim(W) $. Moreover, since $k(\hat{V}) = k(W)[H]$, the degree of $G$ is $[k(\hat{V}):k(W)]$. 
 If $G(T)= \sum a_{i} T^{i}$, we recall that, for $y \in W $, $G(y,T)= \sum a_{i}(y) T^{i}$.

\begin{prop} [Strongly finite morphism] \label{Strongly finite morphism}
Let $ u : V \rightarrow W $ be a strongly finite morphism. Let $\hat{V}$ be given by Proposition \ref{reduction} such that   $k[\hat{V}] = k[W][H]$ where $H $ is  algebraic over $k(W)$ with minimal polynomial $G(T) \in k[W][T]$. Let $n=deg(G)$. Then 
\begin{itemize}
\item[$(a)$] For all $y \in W $,  $\vert u^{-1}(y) \vert = \vert \{ z \in k,  G(y, z)=0  \} \vert  $

\item[$(b)$] for all  $ m \geq 1 $, $ \{ y \in W,    \vert  u^{-1}(y) \vert \geq m \} =  U_{m}(G,W)$

\item[$(c)$] $F_{n}$ is open and for all  $ m \geq 1 $, $F_{m}$ is locally closed.

\end{itemize}
\end{prop}

\paragraph{Proof}

$(c)$ follows from $(b)$ since $(b)$ implies  $F_{m} =  U_{m}(G,W) \setminus U_{m+1}(G,W)$.  
 $(b)$  follows from $(a)$. Therefore, we only need to prove $(a)$.

 Set $n=[k(\hat{V}):k(W)]$ and $G(T)=a_{0}+ a_{1}T +\ldots+a_{n-1} T^{n-1}+T^{n} $, $ a_{j} \in k[W]$.  For $y \in W $,  $ G(y, T) = a_{0}(y)+ a_{1}(y)T +\ldots+a_{n-1}(y) T^{n-1}+T^{n} $. For any $ x \in \hat{V} $, we have $ G(\hat{u}_{1}(x), H(x)) = a_{0}(\hat{u}_{1}(x))+ a_{1}(\hat{u}_{1}(x))H(x) +\ldots+a_{n-1}(\hat{u}_{1}(x)) H^{n-1}(x)+H^{n}(x)$.

Now, we consider $ W \times k $. It is an affine irreducible algebraic variety and  $ dim( W \times k ) = dim( W) +1$ . The ring of regular morphisms of $ W \times k $ is $k[W][T]$. Since $G$ is irreducible, $k[W][T]/G$ is integral and therefore the set $ M = \{ (y,z) \in W \times k, \  G(y, z)=0 \} $ is an irreducible affine algebraic variety. Moreover $dim (M) = dim( W \times k )-1 $ and thus $dim(M) = dim(W) = dim(\hat{V})$.

Since for any $ x \in \hat{V} $,  $G(\hat{u}_{1}(x), H(x))=0 $, it follows that   $ \Psi = (\hat{u}_{1}, H) $ is a regular injective dominant morphism from $\hat{V}$ to $M $.  Clearly $ \Psi$ is finite and therefore it is bijective. 

Let $ \pi : M \rightarrow W $, $ \pi(y,z) = y $, we have $ \hat{u}_{1} = \pi  \circ \Psi $. For any $y \in W $ we have $ \pi^{-1}(y) = \{ (y,z) \in M,  G(y, z)=0  \} $.  Since $ \Psi$ is bijective, we therefore have $ \vert\hat{u}_{1}^{-1}(y) \vert = \vert \pi^{-1}(y)\vert $. Since $ \vert\hat{u}_{1}^{-1}(y) \vert = \vert u^{-1}(y) \vert $
we conclude that, 
for any $y \in W $,  $\vert u^{-1}(y) \vert = \vert \{ z \in k,  G(y, z)=0  \} \vert   $.\\

\subsection{The local  result}

Now, without assuming that the morphism is strongly finite, we will show that we can find an open set on which we have the same properties as in the strongly finite case. This is what we call the local result.

\begin{lemma} [Local result] \label{Local General case}
Let $V$ and $W$ be irreducible  algebraic  varieties over $k$ such that $dim(V) = dim(W) $. Let $ u : V \rightarrow W $ be a dominant morphism. Let $A$ be a non-empty affine open subset of $W$. There exist a non-empty  open subset $U$ of $W$,  $U \subset A$,  and   a monic irreducible   polynomial $G(T) \in k[A][T]$,  such that 

\begin{itemize}
\item[$(a)$] For all $y \in U $,  $\vert u^{-1}(y) \vert = \vert \{ z \in k,  G(y, z)=0  \} \vert  $

\item[$(b)$] for all  $ m \geq 1 $, $ \{ y \in U,  \vert  u^{-1}(y) \vert \geq m \} = U \cap U_{m}(G,A)$

\item[$(c)$]  $U \cap F_{m}$ is locally closed.

\end{itemize}

\end{lemma}

\paragraph{Proof}
Here again, it suffices to show $(a)$.
Let us first consider the affine case. We assume that $V$ and $W$ are affine and that $A=W$. 

Let $\hat{V}$ be given by Proposition \ref{reduction}. Since $k(\hat{V})$ is separable, we can choose $L$ in $k[\hat{V}]$ such that $k(\hat{V})=k(W)[L]$. There is $c \in k[W]$ such that the minimal polynomial $G$ of $H=cL$ is in $k[W][T] $,  $k(\hat{V})=k(W)[H]$ and  $deg(G)=[k(V):k(W)]_{s}$.

 If $k[\hat{V}]= k[W][H]$, the morphism $u$ is strongly finite. We take $U=W$ and the result is given in Proposition~\ref{Strongly finite morphism}. 
 
 Otherwise   there is $s\geq 1$ and $a_{1},\ldots,a_{s}$  such that $k[\hat{V}]= k[W][H][a_{1},\ldots,a_{s}]$. 

Let $W_{1}$ be the irreducible affine algebraic variety associated with the domain $k[W][H]$.  To the injections $k[W] \hookrightarrow k[W][H] \hookrightarrow k[\hat{V}] $, correspond morphisms $ v_{1} :W_{1} \rightarrow W$ and $ v_{0} : \hat{V} \rightarrow W_{1} $ and we have $ \hat{u}_{1} = v_{1} \circ v_{0} $. The morphism $v_{1}$ is separable and strongly finite. The morphism $ v_{0}$ is dominant and we set $Z_{0}$ the closure of the set $ \{ x' \in W_{1}, x' \notin v_{0}(\hat{V}) \} $. Since $v_{0}(\hat{V})$ contains a non-empty open subset of $W_{1}$, $Z_{0}$ is a proper closed subset of $W_{1}$. 

Moreover, since $k(\hat{V})=k(W)[H]$, $ a_{1},\ldots,a_{s} $ are in  $k(W)[H]$. Therefore for any $x \in \hat{V}$, $ a_{i}(x)= \frac{b_{i}(v_{0}(x))}{c_{i}(v_{0}(x))}$.
The closed set $Z_{i} = \{ x' \in W_{1}, c_{i}(x')=0 \}$ is a proper closed subset of $W_{1}$. Let $Z = Z_{0} \cup Z_{1} \cup\ldots\cup Z_{s} $ and $U'= W_{1} \setminus Z $. The set $U'$ is a non-empty open subset of $W_{1}$. Using the injection $ \psi =(v_{0},  a_{1},\ldots,a_{s}) $, it is easy to show that for any $ x' \in U'$, $\vert v_{0}^{-1} (x') \vert =1 $. 
Now, since $v_{1}$ is finite, $v_{1}(Z)$ is a proper closed subset of $W$. We set $U = W \setminus v_{1}(Z)$. We have $v_{1}^{-1} (U) \subset U'$. Therefore, for any $y \in U $, $ \vert u^{-1}(y) \vert = \vert\hat{u}_{1}^{-1}(y) \vert  = \vert v_{1}^{-1} (y) \vert $ and the result in the affine case follows from Proposition~\ref{Strongly finite morphism} applied to $v_{1}$.

Let us now return to the general case. The algebraic varieties $V$ and $W$ are irreducible but are no longer assumed to be affine. 

The set $A$ is a non-empty affine open subset of $W$, therefore it is a non-empty irreducible affine algebraic variety and $dim(A) =dim(W)$. 
Since $u$ is dominant, the set $B = u^{-1}(A) $ is a non-empty open subset of $V$. There is a non-empty affine open subset $A' \subset B$. We denote $u'$ the restriction of $u$ to $A'$. Clearly $u': A' \rightarrow A$ is a dominant morphism and $dim(A') = dim(A)$. Applying the result in the affine case, there exist  a non-empty  open set $U_{1} \subset A $ and   a monic irreducible polynomial $G(T) \in k[A][T]$ with degree $n = [k(A') : k(A)]_{s}$  such that 
for all $y \in U_{1} $,  $\vert u'^{-1}(y) \vert = \vert \{ z \in k,  G(y, z)=0  \} \vert$.

Now let $ C = V \setminus A'$ and $F$ the closure of $u(C)$. Since $dim(C) < dim(V)$ we have $dim(F) < dim(W)$ and therefore its complement $U_{2}$ is a non-empty open subset of $W$. Set $U= U_{1} \cap U_{2}$, then for any $y \in U$, $u^{-1}(y) =  u'^{-1}(y)$ and the result follows.

\subsection{The result on closed subsets.}

Let $V$ and $W$ be algebraic varieties not assumed irreducible and let  $ u : V \rightarrow W $ be a  regular dominant morphism.   

Let $Y$ be a closed subset of $V$ and let $M$ be the closure of $u(Y)$. We say that $Y$ is full if $ u^{-1}(M)=Y $. Note that $Y$ is full if and only if there exists a closed subset $Z$ of $W$ such that $Y=u^{-1}(Z)$. 

We then have the following lemma.

\begin{lemma} [Closed subsets] \label{Constructibility on closed subsets}
Let $Y$ be a closed subset of $V$ and let $M$ the closure of $u(Y)$. Let us assume that  $Y$ is full. Then

\begin{itemize}
\item[$(1)$] there is a non-empty open subset $U$ of $M$ , $U \subset u(V)$,  such that for all $m > 0 $  or $m = \infty$, $F_{m} \cap U$ is a constructible set

\item[$(2)$] there is  an integer $l$, such that  $F_{m} \cap U$ is empty for all $l < m < \infty$

\item[$(3)$]  Let $Y'=  Y \setminus \mathcal{U} $ where $ \mathcal{U} = u^{-1}(U) $. The closed set $Y'$ is full and is a proper subset of $Y$. 
  
\end{itemize}

\end{lemma}

\paragraph{Proof}

We first prove $(3)$. Let $M'$ be the closure of $u(Y')$. We have $ Y' \subset u^{-1}(M')$. Now since $ u^{-1}(M)=Y $ and $ u^{-1}(U) = \mathcal{U}$
we have $u^{-1}(M \setminus U) = Y' $. Therefore $ M' \subset (M \setminus U)$ and $u^{-1}(M') \subset Y' $.  Consequently, $u^{-1}(M') = Y' $ and $Y'$ is full.  Moreover  it is clear that $\mathcal{U}$ is  a non-empty open subset of $Y$ and therefore $Y'$ is a proper subset of $Y$.
We note that $U = u(\mathcal{U}) $.

Let us now show (1) and (2).

Let $X_{1}$,\ldots,$X_{l}$, $l \geq 1$, be the irreducible components of $Y$, let $Z_{i}$ be the closure of $u(X_{i})$ and let $u_{i}$ the restriction of $u$ to $X_{i}$. Note that $Z_{i}$ is irreducible and $dim(X_{i}) \geq dim(Z_{i})$.  Let $m$ such that  $ dim(Z_{m}) = \max_{1\leq i\leq l}(dim(Z_{i}))$. For each $i$, we either have $Z_{i} = Z_{m}$  or $Z_{i} \cap Z_{m}$ is a proper closed subset of $Z_{m}$. 
For convenience, the sets are reindexed so that for $1 \leq i \leq I$ (resp. $ I+1 \leq j \leq I+J$)  $Z_{i} = Z_{m}$ (resp. $Z_{j} \cap Z_{m} \neq Z_{m}$); for $J=0$, there is no set $Z_{j}$.  If  $J \neq 0 $, for each $I+1 \leq j \leq I+J$ let $U_{jm}= (W \setminus Z_{j}) \cap Z_{m}$. The sets $U_{jm}$ are non-empty open subsets of $Z_{m}$. Their intersection  $ G_{m} = \cap_{j \in J} U_{jm}$ is  a non empty open subset of $Z_{m}$ and of $M$. More precisely, $ G_{m} = W \setminus (\cup_{j \in J} Z_{j})$.  When $J=0$, $ G_{m} = Z_{m}$.

There are  two cases.

\paragraph{\textit{Case 1}} There is $i$, $1 \leq i \leq I $, such that $dim(X_{i}) > dim(Z_{m}) $.  In this case, there is a non-empty open subset $U_{\infty} $ of $Z_{m}$, such that for any $y \in U_{\infty}$ we have $ \vert u_{i}^{-1}(y) \vert = \infty $, that is to say $U_{\infty} \subset F_{\infty} $. The set $U = U_{\infty} \cap G_{m} $ is a non-empty open subset of $M$. For $m = \infty $, (1) is satisfied since $U \cap F_{\infty} = U $ is constructible. For $m < \infty $, $U \cap F_{m}$ is empty and therefore constructible. (2) is obvious since $U \cap F_{m}$ is empty for $m < \infty $. 
 
\paragraph{\textit{Case 2}} For all $1 \leq i \leq I $ we have $dim(X_{i}) = dim(Z_{m}) $. For each $i \neq i'$,   $1 \leq i,i' \leq I$, $ dim(X_{i} \cap X_{i'}) < dim(Z_{m}) $. Therefore the  closure in $ Z_{m} $ of $u(X_{i} \cap X_{i'}) $ is a proper subset of $ Z_{m} $ and its complement $D_{ii'}$ is a non-empty  open subset of $Z_{m}$. We denote by $D_{m}$ the intersection of all the $D_{ii'}$ with $i \neq  i'$.  

 For each $ 1 \leq i \leq I $ and each $k \geq 1 $, let us set $F_{ik}= \{ y \in Z_{m},   \vert u_{i}^{-1}(y) \vert = k \}$. Using Lemma~\ref{Local General case} for each $i$, we find an  open subset $U_{i}$ of $Z_{m}$ such that $U_{i} \cap F_{0}$ and $U_{i} \cap F_{\infty}$ are empty and such that  $F_{ik} \cap U_{i}$ is a constructible set that is empty for $k$ greater than some $k_{i}$. We denote by $L_{m}$  the intersection of the $U_{i}$, $1 \leq i \leq I $. Now the set $U = L_{m} \cap G_{m} \cap D_{m} $ is a non-empty  open subset of $M$ and of $Z_{m}$. Of course, $F_{ik} \cap U$ is still a constructible set. Moreover it is clear that, for any $y \in U $, 
 \begin{itemize}
\item[$(a)$]
 for all $I+1 \leq j \leq I+J $, $ u^{-1}(y) \cap X_{j} = \emptyset $
 \item[$(b)$]
  for all $1 \leq i \leq I $, $ u^{-1}(y) \cap X_{i} \neq \emptyset $
 \item[$(c)$]
 for all $i \neq i'$,  $1 \leq i,i' \leq I$, $ u_{i}^{-1}(y) \cap u_{i'}^{-1}(y) = \emptyset $.
 \end{itemize}
  It follows from $(c)$,  that for any $y \in U $, $ u^{-1}(y) $ is the disjoint union of  the $ u_{i}^{-1}(y) $, $1 \leq i \leq I $. 
 
 We note that, in view of $(b)$, $ \vert u^{-1}(y)\vert  \geq I $ ; this means that $F_{k} \cap U $ is empty for $k \leq I$.    
 
 For $ k \geq I $, we have
 $$F_{k} \cap U = \cup  ((F_{1l_{1}} \cap U ) \cap \ldots\cap (F_{Il_{I}} \cap U))$$
 where the union is over the $ (l_{1},\ldots,l_{I} ) $ such that $ \sum_{1 \leq i \leq I}l_{i}=k $.

 We conclude that $ F_{k} \cap U$  is constructible for all $k$, and that it is empty  for $ k \geq I \times \max_{1\leq i\leq I} k_{i} $.  
 Therefore the lemma is proved.

\subsection{End of the proof of Theorem~\ref{General Case}}

We conclude now the proof of Theorem~\ref{General Case}.

\paragraph{Proof}

 We define a sequence $ Y_{1},\ldots,Y_{m} $ of full closed subsets of $V$ as follows.  We first set $ Y_{1} =V $. Then, for each  $s$, the open set $U$ given by Lemma~\ref{Constructibility on closed subsets} applied to $Y_{s}$ is denoted by $U_{s}$, and we set $ \mathcal{U}_{s} = u^{-1}(U_{s}) $ and  $Y_{s+1} = Y_{s} \setminus \mathcal{U}_{s}$.   
 Since the sequence of closed sets $ Y_{1},\ldots,Y_{m}$ is decreasing, there is $m$ such that $Y_{m+1}$ is empty. Since, for each  $s$, $Y_{s}= Y_{s+1} \cup \mathcal{U}_{s}$, it follows that  $V = \mathcal{U}_{1} \cup \ldots\cup \mathcal{U}_{m}$.  Therefore for $j > 0 $  or $j = \infty$, we have $ F_{j} = (U_{1} \cap F_{j} ) \cup\ldots\cup (U_{m} \cap F_{j} ) $. 
 
Since each $(U_{s} \cap F_{j} )$ is constructible, $ F_{j}$ is constructible. 

Moreover, in each step $ s=1,\ldots,m $, there is $l_{s}$ such that $ U_{s} \cap F_{j} $ is empty for $ l_{s} < j <  \infty $. Taking $l = \max_{1\leq s\leq m} l_{s} $, we have $F_{j}$  is empty for $ l < j <  \infty $.

\end{document}